\documentclass[12pt]{article}
\usepackage{amsmath}
\usepackage{amsfonts}
\usepackage{amsthm}
\begin{document}

\def\N{\mathbb{N}}
\def\F{\mathbb{F}}
\def\Z{\mathbb{Z}}
\def\R{\mathbb{R}}
\def\Q{\mathbb{Q}}
\def\H{\mathcal{H}}
\def\fn{\F(q)}
\def\h{\H(q)}
\parindent= 3.em \parskip=5pt
\baselineskip=17pt
\centerline{\bf{ ALGEBRAIC CONTINUED FRACTIONS IN $\F_q((T^{-1}))$}}
 \centerline{\bf{AND RECURRENT SEQUENCES IN $\F_q$ }}

\vskip 0.5 cm
\centerline{\large Alain Lasjaunias \footnote{Alain.Lasjaunias@math.u-bordeaux1.fr}}
\centerline{\small\it C.N.R.S.-UMR 5251, Universit\'e
    Bordeaux I, Talence 33405, FRANCE}
\vskip 1 cm
{\bf{Abstract.}}There exists a particular subset of algebraic power series over a finite field
which, for different reasons, can be compared to the subset of quadratic real numbers. The continued
 fraction expansion for these elements, called hyperquadratic, can sometimes be fully explicited.
In this work, which is a continuation of [L1] and [L2], we describe this expansion for a wide family 
of hyperquadratic power series in odd characteristic. This leads to consider interesting recurrent
 sequences in the finite base field when it is not a prime field.
\vskip 0.5 cm
\noindent \emph{Keywords:} Continued fractions, Fields
of power series, Finite fields.
\newline 2000  \emph{Mathematics Subject Classification:} 11J70, 11T55.
\vskip 1 cm

\centerline{\bf{1. Introduction }}
\vskip 0.5 cm
\par Formal power series over a finite field are analogues of real numbers. Like quadratic
 real numbers, for which the continued fraction expansion is well known, certain algebraic 
power series have a continued fraction expansion which can be explicitly described. 
Most of these power series belong to a particular subset of algebraic elements related to the 
existence of the Frobenius isomorphism in these power series fields. The reader may consult [BL] 
for further information on these elements called hyperquadratic. In a recent work [L1] we have 
introduced a family of hyperquadratic elements having a continued fraction expansion with a
regular pattern. This expansion is linked to particular sequences in a finite field. Here we complete
the study of these sequences. It is also worth mentioning that, in an unexpected way, 
the present work sheds a new light on an older one [LR].  
\par We are concerned with power series over a finite field $\F_q$
of odd characteristic $p$. Given a formal indeterminate $T$, we
consider the ring of polynomials $\F_q[T]$ and the field of rational
functions $\F_q(T)$. Then if $\vert T\vert $ is a fixed real number greater than one, we introduce the ultrametric
 absolute value defined on the field $\F_q(T)$ by $\vert P/Q \vert=\vert T\vert
 ^{\deg(P)-\deg(Q)}$. The completion of this field for this absolute
 value is the field of power series in $1/T$ over
$\F_q$, which is often denoted by $\F_q((T^{-1}))$ or here simply by $\F(q)$.
 If $\alpha \in \F(q)$ and $\alpha \neq 0$, we have
$$\alpha=\sum_{k\leq k_0}u_kT^k,\quad  \text{ where }k_0 \in \Z,
u_k\in \F_q, u_{k_0}\neq 0 \quad \text{ and }\quad 
\vert \alpha\vert=\vert T\vert^{k_0}.$$
We know that each irrational element $\alpha$ of
$\F(q)$ can be expanded as an infinite continued fraction. This will be denoted by
$\alpha=[a_1,\dots,a_n,\dots]$, where the $a_i\in \F_q[T]$ are non
constant polynomials (except possibly for the first one) and are
called the partial quotients of $\alpha$. As usual the tail of the expansion,
$[a_n,a_{n+1},\dots]$, called the complete quotient,
 is denoted by $\alpha_n$ where $\alpha_1=\alpha$. The numerator and the denominator of the
truncated expansion $[a_1,a_2,\dots,a_n]$, which is called a convergent,
are denoted by $x_n$ and $y_n$. These polynomials, called continuants, are
 both defined by the same recursive relation
 :$K_n=a_nK_{n-1}+K_{n-2}$ for $n\geq 2$, with the initials
 conditions $x_0=1$ and $x_1=a_1$ for the the sequence of numerators,
 while the initial conditions are $y_0=0$ and $y_1=1$ for the sequence
 of denominators. For a general account on continued fractions in power
 series fields and also for numerous references the reader may consult
 W. Schmidt's article [S].
\par In this note we consider continued fraction expansions for algebraic
power series over a finite field. We recall that the first works in
this area are due to L. Baum and M. Sweet [BS] and later to W. Mills
and D. Robbins [MR]. The problem discussed here has been introduced in
[L1]. In the next section we present the background of this problem and
we state a technical lemma to go from characteristic zero to positive
characteristic. In the third section, with
Proposition A, we define a large class of algebraic continued fractions in the fields $\F(q)$.
 In the fourth section we state the main result, Theorem B, which gives an
explicit description of these continued fractions under certain conditions. We also present
 an illustration in the field of power series over $\F_{27}$ which is
 Corollary B. At the end of this section we state a conjecture
 concerning a family of irreducible polynomials over $\F_p$. The last section is dedicated to the proof of Theorem B and its corollary. 
\vskip 0.5 cm
\centerline{\bf{2. A special pair of polynomials }}
\vskip 0.5 cm
\par For each integer $k\geq 1$, we consider the following pair of
polynomials in $\Q[T]$:
$$P_k(T)=(T^2-1)^k \quad \text{ and }\quad
Q_k(T)=\int_{0}^{T}(x^2-1)^{k-1}dx.$$ We have the following finite continued
fraction expansions in $\Q(T)$ :
$$P_1(T)/Q_1(T)=[T,-T], \quad P_2(T)/Q_2(T)=[3T,T/3,-3T/4,-4T/3]$$
and more generally
$$P_k/Q_k=[v_{1,k}T,\dots,v_{i,k}T,\dots,v_{2k,k}T],\eqno{(1)}$$
where the rational numbers $v_{i,k}$ for $k\geq 1$ and $1\leq i\leq 2k$
are defined by $v_{1,k}=2k-1$ and recursively, for $1\leq i\leq 2k-1$, by
$$v_{i+1,k}v_{i,k}=(2k-2i-1)(2k-2i+1)(i(2k-i))^{-1}.\eqno{(2)}$$
This continued fraction expansion for $P_k/Q_k$ has been established
in [L1]. Moreover we consider the rational numbers
$$\theta_k=(-1)^k2^{-2k}\binom{2k}{k}\quad \text{ and }\quad
\omega_k=-(2k\theta_k)^{-2} \quad \text{ for }\quad 
k\geq 1.\eqno{(3)}$$
These rational numbers were introduced in [L1] in connection with the
pair $(P_k,Q_k)$. Indeed we have $Q_k(1)=-(2k\theta_k)^{-1}$ and also
$$v_{2k+1-i,k}=v_{i,k}\omega_k^{(-1)^{i+1}}.\eqno{(4)}$$

\par  We recall that throughout this note $p$ is an odd prime
number. Our aim is to obtain, by reducing the identity (1) modulo $p$, a similar identity in $\F_p(T)$. 
Clearly the integer $k$ must be well chosen. The easiest way to do so is to assume that $2k<p$ and 
this is what we did in [L1]. Here we shall extend this to other values of $k$. We set $r=p^t$ where $t$ 
is a positive integer. Then we introduce the subset $E(r)$ of integers $k$ such that 
$$k=mp^l+(p^l-1)/2 \quad \text{ for }\quad 1\leq m\leq (p-1)/2 \quad \text{ and }\quad 0\leq l\leq t-1.\eqno{(5)}$$
For instance $E(3)=\lbrace 1 \rbrace$, $E(5)=\lbrace 1,2 \rbrace$ and
$E(25)=\lbrace 1,2,7,12 \rbrace$. Note that we have $E(r)\subset
\lbrace 1,\dots,(r-1)/2 \rbrace$ with equality if $r=p$. Also
$(r-1)/2\in E(r)$ in all cases. We have the following result where $\bold{v}_p(x)$ is used to denote the p-adic valuation 
of a rational number $x$.
\vskip 0.2 cm
\noindent {\bf{Lemma 1.}}{\emph{ Let $p$ and
    $r$ be as above. Let $k$ be a positive integer with $k\in E(r)$.
 \newline 1) For $1\leq i\leq 2k-1$ we have $\bold{v}_p(i)=\bold{v}_p(2k-2i+1)$ and 
$\bold{v}_p(2k-i)=\bold{v}_p(2k-2i-1)$. For $1\leq i\leq 2k$ we have $\bold{v}_p(v_{i,k})=0$.
 Consequently,in the sequel, for $1\leq i\leq 2k$, $v_{i,k}$  and for
 $1\leq i\leq 2k-1$, $i/(2k-2i+1)$ and $(2k-i)/(2k-2i-1)$ will be considered as elements of $\F_p^*$.
\newline 2) For $0\leq i\leq 2k$ we have
$\bold{v}_p(\binom{2k}{i})=0$. Consequently $\theta_k$, $2k\theta_k$
and $\omega_k$ as well as $\binom{2k}{i}$ for $0\leq i\leq 2k$ will be considered in the sequel
 as elements of $\F_p^*$.
\newline 3) We can define in $\F_p[T]$ the pair of polynomials 
$$P_k(T)=(T^2-1)^k \quad \text{ and }\quad Q_k(T)=\sum_{0\leq i\leq
  k-1}b_iT^{2i+1}$$
where $b_i=(-1)^{k-1-i}\binom{k-1}{i}(2i+1)^{-1} \in \F_p$. The identity (1) holds in $\F_p(T)$
for the rational function $P_k/Q_k$, with the $v_{i,k}$ defined in $\F_p^*$ as above.
}}
\newline Proof : Let $k\in E(r)$. According to (5) we have $2k+1=(2m+1)p^l$ with $3\leq 2m+1\leq p$
and $0\leq l\leq t-1$. For $1\leq i\leq 2k-1$ we have $i<p^t$ and therefore 
$\bold{v}_p(i)\leq \bold{v}_p(2k+1)$. This implies clearly that
 $\bold{v}_p(i)=\bold{v}_p(2k-2i+1)$. We also have $2k-2i-1=2(2k-i)-(2k+1)$ and consequently, 
changing $i$ into $2k-i$, the same arguments show that $\bold{v}_p(2k-i)=\bold{v}_p(2k-2i-1)$. 
By (2), it follows that $\bold{v}_p(v_{i,k}v_{i+1,k})=0$ for $1\leq i\leq 2k-1$. Since
 $\bold{v}_p(v_{1,k})=\bold{v}_p(2k-1)=0$ we have $\bold{v}_p(v_{i,k})=0$ for $1\leq i\leq 2k$. 
So we have proved the first point. For the second one we use a classical formula on the p-adic
valuation of $n!$. Indeed for an integer $n\geq 1$ and a prime $p$ we have 
$\bold{v}_p(n!)=(n-s_p(n))/(p-1)$, where $s_p(n)$ denotes the sum of the digits of $n$ when it 
is written in basis $p$. Since $k\in E(r)$ we can write $2k=2mp^l+(p-1)(p^{l-1}+\dots+1)$. 
For $0\leq i\leq 2k$, this writing implies the equality $s_p(2k-i)+s_p(i)=s_p(2k)$. Consequently we have 
$\bold{v}_p((2k-i)!)+\bold{v}_p(i!)=\bold{v}_p((2k)!)$ and therefore $\bold{v}_p(\binom{2k}{i})=0$
 for $0\leq i\leq 2k$. Now we prove the last point. According to (1),
 by a trivial integration, we can write in $\Q(T)$
$$\sum_{0\leq i\leq k-1}b_iT^{2i+1}=(T^2-1)^k[0,v_{1,k}T,\dots,v_{2k,k}T]$$ 
Since $k\in E(r)$, the right hand side of this equality can be reduced modulo $p$ in $\F_p(T)$. 
Thus the left hand side is well defined by reduction modulo $p$, i.e. $\bold{v}_p(b_i)\geq 0$ for 
$0\leq i \leq k-1$. Consequently the pair $(P_k,Q_k)$ is well defined in $(\F_p[T])^2$ and we have 
the desired continued fraction expansion  for the rational function $P_k(T)/Q_k(T)$ in $\F_p(T)$. 
This completes the proof of the lemma.
\par Given $r$ and $k\in E(r)$, we need now to introduce a pair
of finite sequences $(g_i)_{0\leq i\leq 2k}$ and $(h_i)_{0\leq i\leq
  2k}$ of functions in $\F_p(X)$ which will be used further on.
 We set
$$\left  \{ \aligned & g_0(X)=\theta_k + X,\qquad
  g_{2k}(X)=1/(\theta_k-X) \\& and \quad for \quad  1\leq i\leq 2k-1 \\&  
  g_i(X)=2k\theta_kv_{i,k}(i/(2k-2i+1)) \frac
  {\theta_k+w_{i,k}X}{\theta_k+w_{i-1,k}X} \endaligned \right.\leqno{(G)}$$
also
$$\left  \{ \aligned & h_0(X)=X/(\theta_k + X), \qquad h_{2k}(X)=X/(\theta_k -X) 
\\& and \quad for \quad  1\leq i\leq 2k-1 \\&  
  h_i(X)=(-1)^i\binom{2k}{i}\frac {\theta_k X}{(\theta_k+w_{i,k}X)(\theta_k+w_{i-1,k}X)}
 \endaligned \right.\leqno{(H)}$$
where
$$w_{i,k}=(-1)^i\binom{2k-1}{i}\in \F_p \quad \text{ for }\quad 0\leq i\leq
2k-1.$$
Due to Lemma 1, the functions defined above are not zero. Moreover we may have
$w_{i,k}=0$ for some $i$ but $w_{i,k}-w_{i-1,k}=(-1)^i\binom
{2k}{i}\neq 0$ for $1\leq i\leq 2k-1$.  
\vskip 0.5 cm
\centerline{\bf{3. Continued fractions of type $(r,l,k)$ in $\F(q)$ }}
\vskip 0.5 cm

\par In [L1] a process to generate in $\F(q)$ algebraic
continued fractions from certain polynomials in $\F_q[T]$ is
presented. The following proposition is a
particular case of a more general theorem (see [L1] Theorem 1, p. 332-333).
\vskip 0.2 cm
\noindent {\bf{Proposition A.}}{\emph{ Let $p$ be an odd prime number. Set
    $q=p^s$ and $r=p^t$ with $s,t\geq 1$. Let $k$ be an integer with
    $k\in E(r)$. Let $(P_k,Q_k)\in (\F_p[T])^2$ be defined as in Lemma
    1. Let $l\geq 1$ be an integer. Let
  $(\lambda_1,\lambda_2,\dots,\lambda_l)$ be a $l$-tuple in
  $(\F_q^*)^l$. Let $(\epsilon_1,\epsilon_2)\in (\F_q^*)^2$. There
  exists a unique infinite continued fraction $\alpha=[\lambda_1T,\dots,\lambda_lT,\alpha_{l+1}]\in
 \F(q)$ defined by  
  $$\alpha^r=\epsilon_1 P_k\alpha_{l+1}+\epsilon_2Q_k.$$
This element $\alpha$ is the unique root in $\F(q)$ with $\vert \alpha
\vert \geq \vert T\vert$ of the algebraic equation
$$\qquad y_lX^{r+1}-x_lX^r+(\epsilon_1P_ky_{l-1}-\epsilon_2Q_ky_l)X-
\epsilon_1P_kx_{l-1}+\epsilon_2Q_kx_l=0$$
where $x_l$,$x_{l-1}$,$y_l$ and $y_{l-1}$ are the continuants defined
in the introduction.}}
\par Note that these continued fractions satisfy an algebraic equation
 of a particular type. The reader may consult the introduction of [BL] for a presentation 
of these particular algebraic power series which are called
hyperquadratic. A continued fraction defined as in Proposition A is generated by the pair $(P_k,Q_k)$ for $k\in E(r)$. Such a
continued fraction will be called an expansion of type $(r,l,k)$. 
When the pair $(P_k,Q_k)$ is fixed, this expansion depends on the $l$-tuple
  $(\lambda_1,\lambda_2,\dots,\lambda_l)$ in $(\F_q^*)^l$ and on the pair
  $(\epsilon_1,\epsilon_2)$ in $(\F_q^*)^2$. When these $l+2$ elements
  in $\F_q^*$ are taken arbitrarily then the expansion has a regular
  pattern only up to a certain point (see [L1] Proposition 4.6, p. 347). In
  the next section we are concerned with a particular subfamily of
  these continued fractions. 
\vskip 0.5 cm
\centerline{\bf{4. Perfect continued fractions of type $(r,l,k)$ in $\F(q)$ }}
\vskip 0.5 cm
In previous works we have seen that an expansion of type
$(r,l,k)$, under certain conditions on
$(\lambda_1,\lambda_2,\dots,\lambda_l)$ and $(\epsilon_1,\epsilon_2)$,
may be given explicitly. A first example was given in [L1] Theorem
3. In [L2] a more general case was treated, but there we restricted
ourselves to the case of a prime base field $\F_p$ and we also only
considered the case $r=p$. Note that in this way we could prove the conjecture for the
expansion of a quartic power series over $\F_{13}$ made by Mills and Robbins in [MR] p. 403.
Here our aim is to describe explicitely many expansions of type $(r,l,k)$ having a very regular
 pattern as Mills and Robbins' example does. To do so we need first to introduce
  further notations. Given $l\geq1$ and $k\geq 1$, we define the
  sequence of integers $(f(n))_{n\geq 1}$ where
  $f(n)=(2k+1)n+l-2k$. We also define the sequence of integers $(i(n))_{n\geq 1}$ in the
following way :
$$i(n)=1 \quad \text{if }\quad n\notin f(\N^*)\quad \text{and }\quad i(f(n))=i(n)+1.$$
Finally we introduce the sequence $(A_i)_{i\geq 1}$ of polynomials in
$\F_p[T]$ defined recursively by
 $$A_1=T \quad \text{ and } \quad A_{i+1}=[A_i^r/P_k] \quad \text{
   for }\quad i\geq 1$$
(here the square brackets denote the integer part, i.e. the polynomial
part). Note that the sequence $(A_i)_{i\geq 1}$ depends on the polynomial $P_k$ chosen 
with $k\in E(r)$. It is remarkable that if $2k=r-1$ then this sequence of polynomials is constant
 and we have $A_i=T$ for $i\geq 1$. 
\par For an arbitrary continued fraction of type $(r,l,k)$ the
sequence of partial quotients is based on the above sequence
$(A_i)_{i\geq 1}$ but only up to a certain rank (see the remark after
Lemma 5.1 below). Nevertheless it may happen that this sequence of
partial quotients is entirely described by means of this sequence
$(A_i)_{i\geq 1}$. The aim of the following theorem is to give this
description as well as the conditions of its existence. These
particular expansions of type $(r,l,k)$, which are defined in this
 theorem, will be called perfect (This term was introduced in [L1], p 348).
\vskip 0.2 cm
\noindent {\bf{Theorem B.}}{\emph{ Let $p$ be an odd prime and $q=p^s$, $r=p^t$ with $s,t\geq 1$ 
be given. Let $k\in E(r)$. Let $(A_i)_{i\geq 1}$ in $\F_p[T]$, $(f(n))_{n\geq 1}$ and
$(i(n))_{n\geq 1}$ in $\N^*$ be the sequences defined above. Let $\alpha \in \F(q)$ be a continued fraction of type $(r,l,k)$
 defined by the l-tuple $(\lambda_1,\dots,\lambda_l)$ in $(\F_q^*)^l$ and by the pair
  $(\epsilon_1,\epsilon_2)$ in $(\F_q^*)^2$.
 Then the partial quotients of this expansion satisfy 
$$(I)\quad a_n=\lambda_nA_{i(n)}\quad \text{ where }\quad \lambda_n \in
\F_q^* \quad \text{ for }\quad n\geq 1$$
 if and only if we can define in $\F_q^*$
$$(II)\qquad  \delta_n=2k\theta_k [\lambda_n^r,\dots,\lambda_1^r,2k\theta_k\epsilon_2^{-1}] \quad
\text{ for }\quad 1\leq n\leq l$$ 
and we have $(III)$  
$$\text{ either case }\quad (III_1):\qquad \delta_l=4k^2\theta_k(\epsilon_1/\epsilon_2)^r$$
$$\text{ or case }\quad(III_2):\qquad \delta_l\neq
4k^2\theta_k(\epsilon_1/\epsilon_2)^r$$ and there exists in $\F_q^*$ a sequence
$(\gamma_n)_{n\geq 1}$ defined recursively by 
$$\left  \{ \aligned &
  \gamma_1^r=(4k^2\theta_k\epsilon_1^r\delta_l^{-1}-\epsilon_2^r)\theta_k\delta_1^{-r}  \\& \gamma_n=\gamma_{n-1}(\delta_n\delta_{n-1}\omega_k)^{-1} \quad for \quad 2\leq n\leq l \\& \gamma_{f(n)+i}=C_0h_i(\gamma_n^{r}) \quad
 for \quad 0\leq i\leq 2k\quad  and \quad n\geq 1 \endaligned
\right.\leqno{(\Gamma)}$$
where
$$C_0=\gamma_l\epsilon_1^r(\delta_1\gamma_1)^{-r}(\delta_l\omega_k)^{-1}
\in \F_q^*.$$
\newline If $(II)$ and $(III)$ hold then we can define recursively a sequence
 $(\delta_n)_{n\geq 1}$ in $\F_q^*$ by the
initial values $\delta_1,\dots,\delta_l$ given by $(II)$ and the formulas
 $$\delta_{f(n)+i}=\epsilon_1^{r(-1)^{n+i}}\delta_n^{r(-1)^i}g_{i,n}\quad
\text{for}\quad  n\geq 1\quad
\text{and for}\quad  0\leq i\leq 2k,\eqno{(D)}$$
where $g_{i,n}=g_i(0)$ in case $(III_1)$ and $g_{i,n}=g_i(\gamma_n^r)$ in case $(III_2)$. 
Then the sequence $(\lambda_n)_{n\geq 1}$ in
$\F_q^*$, introduced in $(I)$, is defined recursively by the first values $\lambda_1,\lambda_2,\dots,\lambda_l$ and 
the formulas 
$$\lambda_{f(n)}=\epsilon_1^{(-1)^n} \lambda_n ^r,\quad \lambda_{f(n)+i}
=-v_{i,k}\epsilon_1^{(-1)^{n+i}}\delta_n^{(-1)^i} \eqno{(LD)}$$
for $n\geq 1$ and for $1\leq i\leq 2k$.}}
\par In this theorem we have two conditions $(II)$ and $(III)$ which are not at the same
level. Condition $(II)$ is primary and clearly necessary to define recursively the sequence
$(\delta_n)_{n\geq 1}$ in $\F_q^*$ by $(D)$. This condition has already been pointed out in
[L2] even though there we had only considered the simplest case where the
base field is prime, that is $q=p$. Here it is necessary to underline that there is a
 mistake in the formula given there for $\delta_n$ when $1\leq n\leq l$. Indeed in [L2] Theorem 1 condition
$(H_1)$, instead of $\delta_i=[2k\theta_k\lambda_i,\dots,2k\theta_k\lambda_1,\epsilon_2^{-1}]$ 
one should read  $\delta_i=2k\theta_k[\lambda_i,\dots,\lambda_1,2k\theta_k\epsilon_2^{-1}]$.
 Note that this last formula is in agreement with $(II)$ in Theorem B, if the base field is prime 
and consequently the Frobenius isomorphism is reduced to the identity in $\F_p$. We must also
add that this mistake has no consequence on Theorem 2 of [L2]
because there the values of $\delta_i$ were computed with the right
formula. Condition $(III)$ is of a different kind. It is split into two
distinct cases. In fact case $(III_1)$ has also already
been considered in [L2] Theorem 1 condition $(H_2)$, again when the
base field is prime. 
\par We want now to discuss about case $(III_2)$. This one is more
complex because of a possible obstruction in the recursive definition of the sequence
$(\gamma_n)_{n\geq 1}$ in $\F_q^*$. Actually it is
conjectured that this second case can only happen if
the base field $\F_q$ is a particular algebraic extension of the prime field $\F_p$.
 Indeed the sequence $(\gamma_n)_{n\geq 1}$ is clearly well defined
 if $\gamma_n$ does not belong to $\F_p$ for all $n$. There is a sufficient condition 
to obtain that. We recall that the functions $h_i$ for $0\leq i\leq 2k$ involved in 
the recursive definition of the sequence $(\gamma_n)_{n\geq 1}$ are of two types :
 $h(x)=ax/(x+u)$ or $h'(x)=ax/(x+u_1)(x+u_2)$
where $a\in \F_p^*$ and $u,u_1,u_2\in \F_p$. Consequently if $x$ is an algebraic element 
over $\F_p$ of degree $d>2$ then $h_i(x)$ has degree $d$ (for all $h_i$ of type $h$)
or $d/2$ (eventually for some $h_i$ of type $h'$). This remark implies that if the first $l$
 terms of the sequence $\gamma_n$ have each one a degree over $\F_p$ different 
from a power of two and the constant $C_0$ has a degree over $\F_p$ which is a power of two,
then by induction the degrees over $\F_p$ of all the terms remain greater than one and thus 
none of these terms belongs to $\F_p$. This condition implying the existence of 
the sequence $(\gamma_n)_{n\geq 1}$ may also be necessary but this remains a conjecture.
 Observe that if this conjecture is true and if the base field is prime then the continued 
fraction can only be perfect in case $(III_1)$. We will make a more
precise conjecture in that direction at the end
 of this section. Before going further on, we need to point out the similarity with the problem
 discussed in [LR], particularly on pages 562-565. In this older work we had investigated the
 existence of algebraic continued fractions having linear partial
 quotients, and this matches the case $2k=r-1$ in the present work. The approach in [LR] was singular and completely different
 from here, this forced us to make the restriction $l\geq r$.
\par Now we want to illustrate the occurence of case $(III_2)$ if the base field is $\F_q$ where
 $q=p^m$ and $m$ is not a power of two. We take $p=3$, $q=27$ and
 $r=3$, with $l=1$ and  $k=1$. Since $2k=r-1$, if the expansion is
 perfect, then all partial quotients are linear. The elements of the finite field $\F_{27}$ will be
represented by means of a root $u$ of the irreducible polynomial over
$\F_3$: $P(X)=X^3+X^2-X+1$. Then we have $u^{13}=-1$ and
$$\F_{27}=\lbrace 0,\pm u^i, \qquad 0\leq i\leq 12 \rbrace .$$ 
We have the following corollary.
\vskip 0.2 cm
\noindent {\bf{Corollary C. }}{\emph{ Define the sequences $(\gamma_n)_{n\geq 1}$ and
    $(\delta_n)_{n\geq 1}$ in $\F_{27}^*$ as follows. The
    first is defined recursively by $\gamma_1=u$ and
$$\gamma_{3n-1}=\frac{\gamma_n^3}{1+\gamma_n^3}, \quad
 \gamma_{3n}=\frac{\gamma_n^3}{1-\gamma_n^6},
\quad \gamma_{3n+1}=\frac{\gamma_n^3}{1-\gamma_n^3}\quad
\text{for}\quad n\geq 1.$$
The second one, based upon the first one, is defined recursively by $\delta_1=u^4$
 and
$$\delta_{3n-1}=u^{5(-1)^n}\delta_n^3(1+\gamma_n^3), \quad
 \delta_{3n}=\frac{\gamma_n^3-1}{\delta_{3n-1}},
\quad \delta_{3n+1}=\frac{\delta_{3n-1}}{1-\gamma_n^6}\quad
\text{for}\quad n\geq 1.$$
We consider the following algebraic equation with coefficients in $\F_{27}[T]$ 
$$X^4-TX^3-u^3TX+uT^2-u^6=0.\leqno{(E)}$$
This equation has a unique root $\alpha$ in $\F(27)$ which can be expanded as
an infinite continued fraction
$$\alpha=[T,u^7T,u^2T,u^{11}T,-uT,\dots]=[\lambda_1T,\dots,\lambda_nT,\dots],$$
where the sequence $(\lambda_n)_{n\geq 1}$ in $\F_{27}^*$ is defined recursively by
$\lambda_1=1$ and
$$\lambda_{3n-1}=-u^{6(-1)^n}\lambda_n^3 ,\quad
 \lambda_{3n}=u^{6(-1)^{n+1}}\delta_n^{-1},
\quad \lambda_{3n+1}=-\lambda_{3n}^{-1}\quad
\text{for}\quad n\geq 1.$$}}
 \par Before concluding this section, we make a conjecture in connection with the sequence $(\gamma_n)_{n\geq 1}$ 
described in Theorem B. For the sake of shortness we take $k=1$. Let $p$ be an odd prime. Let us 
consider the three elements of $\F_p(x)$
$$h_0(x)=\frac{2x}{2x-1},\quad h_1(x)=\frac{4x}{1-4x^2} \quad
\text{and }\quad h_2(x)=\frac{-2x}{2x+1}.$$
We define recursively a sequence $(u_n)_{n\geq 1}$ of rational functions in $\F_p(x)$ by
$$u_1(x)=x \quad u_{3n+i-1}(x)=h_i(u_n(x)) \quad \text{for } \quad
0\leq i\leq 2  \quad \text{and } \quad n\geq 1.$$
Let $\mathcal{P}(p)\subset \F_p[x]$ be the subset of all monic polynomials irreducible over $\F_p$
 which appear as a prime factor of the numerator or denominator of $u_n(x)$ for all $n\geq 1$. 
Let $\mathcal{P}_2(p)\subset \F_p[x]$ be the subset of all monic polynomials irreducible over $\F_p$
 of degree $2^k$ for $k\geq 0$. Then the arguments developed after Theorem B show that we 
have $\mathcal{P}(p)\subset \mathcal{P}_2(p)$. We conjecture that
$\mathcal{P}(p)=\mathcal{P}_2(p)$ holds for all odd primes $p$.
   
\vskip 0.5 cm
\centerline{\bf{5. Proofs of Theorem B and Corollary C }}
\vskip 0.5 cm

In this section $p$, $q$ and $r$ are given as above. We consider an integer $k$ with 
$k\in E(r)$ and an integer $l\geq 1$. Moreover the numbers $\theta_k,\omega_k \in
\F_p^*$ and the natural integers $f(n), i(n)$ for $n\geq 1$ are defined as above. The proof of Theorem B will
 be divided into several steps.
\vskip 0.2 cm
\noindent {\bf{Lemma 5.1.}}{\emph{ Let $\alpha=[\lambda_1T,\dots,\lambda_lT,\alpha_{l+1}]\in
 \F(q)$ be a continued fraction of type $(r,l,k)$ for the pair
 $(\epsilon_1,\epsilon_2)\in (\F_q^*)^2$.
 Then there exists a sequence $(\lambda_n)_{n\geq 1}$ in
$\F_q^*$ such that we have 
$$(I)\quad a_n=\lambda_nA_{i(n)}\quad \text{ for }\quad n\geq 1$$
if and only if there exists a sequence $(\delta_n)_{n\geq 0}$ in
$\F_q^*$ such that we have
$$\lambda_{f(n)}=\epsilon_1^{(-1)^n} \lambda_n ^r,\quad \lambda_{f(n)+i}
=-v_{i,k}\epsilon_1^{(-1)^{n+i}}\delta_n^{(-1)^i}\eqno{(LD)}$$ for
$1\leq i\leq 2k$ and $n\geq 1$, with 
$$\delta_n=2k\theta_k^{i(n)}\lambda_n^r-(\omega_k\delta_{n-1})^{-1} \quad
\text{for}\quad  n\geq 1,  \eqno{(D_1)}$$ where $\delta_0=-(\omega_k\epsilon_2)^{-1}$.
}}
\vskip 0.2 cm
\par This lemma, which is the first and the main step in the proof of
Theorem B, is a direct consequence of [L1] Prop 4.6, p 347. There we
proved that an expansion of type $(r,l,k)$ for an arbitrary pair
$(\epsilon_1,\epsilon_2) \in (\F_q^*)^2$ has the pattern given by
$(I)$, where the sequence $(\lambda_n)_{n\geq 1}$ is described by $(D_1)$
and $(LD)$, but only up to a certain rank (if $\delta_n$ ever vanishes
in $(D_1)$). Remark that in the proof of Proposition 4.6 we made the restriction
$2k<p$. This condition was sufficient to have in $\F_p(T)$ the identity $(1)$ of
Section 2
 which is the fundament of the proof.
But, according to Lemma 1 of Section 2, we may replace this condition by $k \in E(r)$ and this
 has no consequences for the proof of the proposition. Now to separate
 the sequence $(\delta_n)_{n\geq   0}$ from the sequence $(\lambda_n)_{n\geq 1}$, we have the following lemma.
\vskip 0.2 cm
\noindent {\bf{Lemma 5.2.}}{\emph{ Let $(\lambda_n)_{n\geq 1}$ and
    $(\delta_n)_{n\geq 0}$ be two sequences in
$\F_q^*$. We assume that they satisfy $(LD)$. Then these sequences
satisfy $(D_1)$ if and only if 
$$(II_0)\qquad
\delta_n=2k\theta_k[\lambda_n^r,\dots,\lambda_{1}^r,\delta_0/(2k\theta_k)]\quad
\text{for}\quad 1\leq n\leq l,$$ 
$$\delta_{f(n)}+(\omega_k\delta_{f(n)-1})^{-1}=\theta_k\epsilon_1^{r(-1)^n}(\delta_n^r+(\omega_k\delta_{n-1})^{-r})\quad
\text{for}\quad  n\geq 1\eqno{(D_2)}$$
and 
$$\delta_{f(n)+i}+(\omega_k\delta_{f(n)+i-1})^{-1}=-2k\theta_kv_{i,k}\epsilon_1^{r(-1)^{n+i}}\delta_n^{r(-1)^i}\eqno{(D_3)}$$
for $1\leq i\leq 2k$ and for $n\geq 1$.}}

\noindent Proof: First we assume that $(D_1)$ holds for $n\geq 1$. For $1\leq n\leq l$ we have $i(n)=1$. By (3) from Section 2, we have
$\omega_k=-(2k\theta_k)^{-2}$, consequently $(D_1)$ for $1\leq n\leq l$ can be written as    
$$\delta_n=2k\theta_k\lambda_n^r+(2k\theta_k)^2\delta_{n-1}^{-1}.$$
By induction, it is clear that $(D_1)$ for $1\leq
n\leq l$ is equivalent to $(II_0)$. Now for $1\leq j\leq 2k$ and for
$n\geq 1$ we have $i(f(n)+j)=1$ consequently, taking into account $(LD)$, we have the equivalence
 between $(D_3)$ and $(D_1)$ at the rank
  $f(n)+j$. Observing
  that for $n\geq 1$ we have $i(f(n))=i(n)+1$ and taking into account
  $(LD)$, we see that $(D_1)$ at the rank $f(n)$ and $n$ implies
  $(D_2)$. Reciprocally suppose $(II_0)$, $(D_2)$ and $(D_3)$
are satisfied. Then $(D_1)$ holds for $1\leq n\leq l$ and also at the rank $f(n)+i$, for $n\geq 1$ and
  $1\leq i\leq 2k$. On the other hand, if $(D_1)$ and $(D_2)$ hold at
  the rank $n\geq 1$, taking into
  account $(LD)$, then $(D_1)$ holds at the rank $f(n)$. Hence, with
  the cases already established and using induction, we see that
  $(D_1)$ holds for $n\geq 1$. This completes the proof of the lemma.
\par In the next lemma we introduce a new sequence $(\gamma_n)_{n\geq 1}$ in $\F_q$
 which is linked to the sequence $(\delta_n)_{n\geq
  0}$.
\vskip 0.2 cm
\noindent {\bf{Lemma 5.3.}}{\emph{ Let $(g_i)_{0\leq i\leq 2k}$ be the
    sequence of functions in $\F_p(X)$ defined by $(G)$ in Section 2. Let $(\delta_n)_{n\geq 0}$ be a sequence in
$\F_q^*$ with $\delta_0,\delta_1,\dots,\delta_l$ given. Then $(\delta_n)_{n\geq 0}$ 
satisfies $(D_2)$ and $(D_3)$ if and only if there exists a sequence
$(\gamma_n)_{n\geq 1}$ in $\F_q$ such that we have
$$\delta_{f(n)+i}=\epsilon_1^{r(-1)^{n+i}}\delta_n^{r(-1)^i}g_i(\gamma_n^r)\quad
\text{for}\quad  0\leq i\leq 2k\quad
\text{and}\quad  n\geq 1,\eqno{(D)}$$
with
$$\gamma_1^r=(\theta_k\delta_0^{-r}-\epsilon_1^r\delta_l^{-1})(\omega_k\delta_1)^{-r}
\eqno{(\Gamma_1)}$$
 and
$$\gamma_n=\gamma_{n-1}(\delta_n\delta_{n-1}\omega_k)^{-1} \quad
\text{for}\quad  n\geq 2.\eqno{(\Gamma_2)}$$
}}
\noindent Proof: First we prove that the sequence $(g_i)_{0\leq i\leq
  2k}$ in $\F_p(X)$, described in $(G)$, can also be defined recursively by $g_0(X)=\theta_k+X$ and 
$$g_{i+1}(X)=2k\theta_k(-v_{i+1,k}+2k\theta_k/g_{i}(X))\quad
\text{for}\quad  0\leq i\leq 2k-1. \eqno{(6)}$$
For $i=0$, $(6)$ becomes
$$g_1(X)=2k\theta_k(-v_{1,k}+2k\theta_k/(\theta_k+X)).$$
Since $v_{1,k}=2k-1$, this equality implies
$$g_1(X)=2k\theta_k(\theta_k-(2k-1)X)/(\theta_k+X).$$
This is in agreement with $(G)$ for $i=1$.
Now we use induction on $i$. Let $1\leq i<2k-1$ and assume that
$g_i(X)$ is as stated in $(G)$. Then, by (2) from Section 2, we have
$$2k\theta_k/g_{i}(X)=\frac{(2k-i)v_{i+1,k}}{2k-2i-1}\frac{\theta_k+w_{i-1,k}X}{\theta_k+w_{i,k}X}.\eqno{(7)}$$
Besides, a direct computation shows that, for $1\leq i<2k-1$, we also have 
$$(2k-i)w_{i-1,k}-(2k-2i-1)w_{i,k}=(i+1)w_{i+1,k}. \eqno{(8)}$$
Combining $(6)$, $(7)$ and $(8)$, we get $g_{i+1}(X)$ as stated in $(G)$.
It remains to compute $g_{2k}(X)$. From $(6)$ and $(7)$ for $i=2k-1$, we obtain
$$g_{2k}(X)=2k\theta_k(-v_{2k,k}-\frac{v_{2k,k}}{2k-1}\frac{\theta_k+(2k-1)X}{\theta_k-X}).\eqno{(9)}$$
Recalling (3) and (4) from Section 2, we also have
$v_{2k,k}=v_{1,k}\omega_k$ and consequently
$$v_{2k,k}=-(2k-1)(2k\theta_k)^{-2}. \eqno{(10)}$$ Finally, combining (9)
and (10), we get  $g_{2k}(X)=1/(\theta_k-X)$ and this is in agreement with $(G)$ for $i=2k$. 
\newline We set now 
$$g_{i,n}=\delta_{f(n)+i}/(\epsilon_1^{r(-1)^{n+i}}\delta_n^{r(-1)^i})\quad
\text{for}\quad  0\leq i\leq 2k\quad
\text{and for}\quad  n\geq 1.\eqno{(11)}$$
Then we define the sequence $(\gamma_n)_{n\geq 1} \in
\F_q$ from the sequence $(\delta_n)_{n\geq 1} \in
\F_q^*$ by
$$\gamma_n^r=g_{0,n}-\theta_k\quad
\text{ for}\quad  n\geq 1. \eqno{(12)}$$
 By definition, $(12)$ becomes $g_{0,n}=g_0(\gamma_n^r)$, thus $(11)$
 implies $(D)$ for $i=0$. Now we prove that $(D_3)$ is equivalent to $(D)$ for
 $1\leq i\leq 2k$. According to
$(11)$, we need to prove that $(D_3)$ is equivalent to $g_{i,n}=g_i(\gamma_n^r)$ for $1\leq i\leq
2k$ and $n\geq1$. Using (11) and again $\omega_k=-(2k\theta_k)^{-2}$, 
$(D_3)$ can be written as  
$$g_{i,n}=2k\theta_k(-v_{i,k}+2k\theta_k/g_{i-1,n}) \quad \text{ for
}\quad 1\leq i\leq 2k. $$
Since $g_{0,n}=g_0(\gamma_n^r)$, with the recursive definition of the
sequence $(g_i)_{1\leq i\leq 2k}$, we see that $(D_3)$ is equivalent to 
$g_{i,n}=g_i(\gamma_n^r)$ for $1\leq i\leq 2k$.
Now we shall see that $(D_2)$ and $(D_3)$ imply $(\Gamma_1)$ and
$(\Gamma_2)$. Hence, with $(\gamma_n)_{n\geq 1}$ defined by $(12)$, we
have $(D)$. For $n=1$, $(D_2)$ becomes
$$\delta_{l+1}+(\omega_k\delta_l)^{-1}=\theta_k\epsilon_1^{-r}(\delta_1^r+(\omega_k\delta_0)^{-r}).\eqno{(13)}$$
But, using $(D)$ for $n=1$ and for $i=0$, we also have  
$$\delta_{l+1}=\epsilon_1^{-r}\delta_1^r(\theta_k+\gamma_1^r).\eqno{(14)}$$
Combining (13) and (14), we obtain the value for $\gamma_1^r$ stated in
$(\Gamma_1)$. Now we assume that $n\geq 2$ and we recall that we have $f(n)-1=f(n-1)+2k$. Consequently,
using $(D)$ for $i=0$ and for $i=2k$, by $(11)$, $(D_2)$ implies  
$$g_{0}(\gamma_n^r)\delta_n^r+(\omega_kg_{2k}(\gamma_{n-1}^r)\delta_{n-1}^r)^{-1}=\theta_k(\delta_n^r+(\omega_k\delta_{n-1})^{-r}).
\eqno{(15)}$$
Since, by $(G)$, $g_{0}(\gamma_n^r)=\theta_k+\gamma_n^r$
and $g_{2k}(\gamma_{n-1}^r)=1/(\theta_k-\gamma_{n-1}^r$), $(15)$ gives 
 $$(\gamma_n\delta_n-\omega_k^{-1}\gamma_{n-1}\delta_{n-1}^{-1})^r=0$$
 which is $(\Gamma_2)$. Reciprocally we assume that both sequences satisfy $(D)$, $(\Gamma_1)$
 and $(\Gamma_2)$. First, as we have seen above, $(D_3)$ hold for $1\leq i \leq
 2k$ and $n\geq 1$. Then $(D)$ for $n=1$ and $i=0$ implies
 $\theta_k+\gamma_1^r=\delta_{l+1}\epsilon_1^r\delta_1^{-r}$. Taking
 $(\Gamma_1)$ into account, this implies $(D_2)$ for $n=1$. Finally, for $n\geq 2$,
  we have seen that $(\Gamma_2)$ implies $(15)$. Using $(D)$ for
  $i=0$ and for $i=2k$, $(15)$ is equivalent to $(D_2)$. The proof of the lemma is complete.
\par In the last lemma we describe the sequence $(\gamma_n)_{n\geq
  1}$, if it is not identically zero.
\vskip 0.2 cm
\noindent {\bf{Lemma 5.4.}}{\emph{ Let $(h_i)_{0\leq i\leq 2k}$ be the
    sequence of functions in $\F_p(X)$ defined by $(H)$ in Section
    2. Let $(\delta_n)_{n\geq 1}$ and 
$(\gamma_n)_{n\geq 1}$ be two sequences in $\F_q^*$ with
$\delta_1,\dots,\delta_l$ and $\gamma_1$ given. We assume that
they satisfy $(D)$. Then they satisfy $(\Gamma_2)$ if and only if we have 
$$\gamma_n=\gamma_{n-1}(\delta_n\delta_{n-1}\omega_k)^{-1} \quad for
\quad 2\leq n\leq l\eqno{(\Gamma'_2)}$$
and $$\gamma_{f(n)+i}=C_0h_i(\gamma_n^{r}) \quad
 for \quad 0\leq i\leq 2k\quad  and\quad for \quad n\geq 1,\eqno{(\Gamma_3)}$$
where
$$C_0=\gamma_l\epsilon_1^r(\delta_1\gamma_1)^{-r}(\delta_l\omega_k)^{-1}
\in \F_q^*.$$
}}
\noindent Proof: First we prove that $(\Gamma_2)$ implies $(\Gamma_3)$. 
We will use the connection between the two sequences
 $(g_i)_{0\leq i\leq 2k}$ and $(h_i)_{0\leq i\leq 2k}$ in $\F_p(X)$. 
Indeed, from $(G)$ and $(H)$, an elementary calculation shows that we have
$$g_i(X)g_{i-1}(X)\omega_k=h_{i-1}(X)/h_i(X)\quad for \quad 1\leq
i\leq 2k.\eqno{(16)}$$
We also have
$$g_0(X)h_0(X)=X \quad \text{and}\quad h_{2k}(X)=Xg_{2k}(X).\eqno{(17)}$$ 
For $1\leq i\leq 2k$ and $n\geq 1$, using $(D)$ and $(16)$, we
have
$$\omega_k\delta_{f(n)+i}\delta_{f(n)+i-1}=h_{i-1}(\gamma_n^r)/h_{i}(\gamma_n^r).\eqno{(18)}$$
Applying $(\Gamma_2)$ at the rank $f(n)+i$, $(18)$ implies
$$\gamma_{f(n)+i}/\gamma_{f(n)+i-1}=h_{i}(\gamma_n^r)/h_{i-1}(\gamma_n^r).\eqno{(19)}$$
Clearly, for $0\leq i\leq 2k$ and $n\geq 1$, from $(19)$ we obtain
$$\gamma_{f(n)+i}=\gamma_{f(n)}h_i(\gamma_n^r)/h_{0}(\gamma_n^r).\eqno{(20)}$$
We assume now that $n\geq 2$. Recalling that $f(n)-1=f(n-1)+2k$, by $(D)$ for $i=0$ and $i=2k$
 and $(17)$, we also have
$$\delta_{f(n)}\delta_{f(n)-1}=(\delta_n\delta_{n-1})^r(\gamma_n/\gamma_{n-1})^r
h_{2k}(\gamma_{n-1}^r)/h_0(\gamma_n^r).\eqno{(21)}$$
Applying $(\Gamma_2)$ at the rank $n$, (21) becomes
$$\omega_k\delta_{f(n)}\delta_{f(n)-1}=h_{2k}(\gamma_{n-1}^r)/h_0(\gamma_n^r).\eqno{(22)}$$
Applying $(\Gamma_2)$ at the rank $f(n)$, $(22)$ implies  
$$\gamma_{f(n)}=\gamma_{f(n)-1}h_{0}(\gamma_n^r)/h_{2k}(\gamma_{n-1}^r).\eqno{(23)}$$
By $(20)$ we also have 
$$\gamma_{f(n)-1}=\gamma_{f(n-1)+2k}=\gamma_{f(n-1)}h_{2k}(\gamma_{n-1}^r)/h_{0}(\gamma_{n-1}^r).\eqno{(24)}$$
Combining $(23)$ and $(24)$, we obtain
$$\gamma_{f(n)}=\gamma_{f(n-1)}h_{0}(\gamma_n^r)/h_{0}(\gamma_{n-1}^r).\eqno{(25)}$$
Consequently, by $(25)$, for $n\geq 1$ we have 
$$\gamma_{f(n)}/h_{0}(\gamma_n^r)=C_0=\gamma_{f(1)}/h_{0}(\gamma_{1}^r).\eqno{(26)}$$
To compute $C_0$, we apply $(\Gamma_2)$ at the rank $f(1)=l+1$ and $(D)$
for $n=1$ and for $i=0$.
 We obtain
$$C_0=\gamma_{l}\epsilon_1^r(\delta_{1}\gamma_1)^{-r}(\delta_{l}\omega_k)^{-1}.$$
Finally, combining $(20)$ and $(26)$, we get $(\Gamma_3)$. We now prove
that $(\Gamma_2')$ and $(\Gamma_3)$ imply $(\Gamma_2)$. Hence
$(\Gamma_2)$ holds for $2\leq n\leq l$. Moreover we obtain $(19)$ directly from
$(\Gamma_3)$. Together with $(18)$, this proves that $(\Gamma_2)$
holds at the rank $f(n)+i$ for $n\geq 1$ and $1\leq i\leq 2k$. We
observe that $(\Gamma_2)$ also holds for $l+1$. Indeed applying
$(\Gamma_3)$ for $n=1$ and $i=0$, together with the value of $C_0$, we
obtain $(\Gamma_2)$ for $l+1$. Now we assume that $n\geq 2$ and we
apply $(\Gamma_3)$ for $i=0$ and for $i=2k$. We have
$$\gamma_{f(n)}/\gamma_{f(n)-1}=h_{0}(\gamma_n^r)/h_{2k}(\gamma_{n-1}^r).\eqno{(27)}$$
Combining $(21)$ and $(27)$ we obtain
$$\gamma_{f(n)}/\gamma_{f(n)-1}=(\delta_n\delta_{n-1})^r(\gamma_n/\gamma_{n-1})^r(\delta_{f(n)}\delta_{f(n)-1})^{-1}.$$
This shows that if $(\Gamma_2)$ holds at the rank $n\geq 2$ then it holds at
the rank $f(n)$. Consequently, with the cases already established and
using induction, we see that $(\Gamma_2)$ holds for
$n\geq 2$. The proof of the lemma is complete.
\vskip 0.2 cm
\noindent {\emph{ Proof of Theorem B:}}
\par Let $\alpha \in \F(q)$ be a continued fraction of type $(r,l,k)$ defined by the $l$-tuple
 $\lambda_1,\dots,\lambda_l \in (\F_q^*)^l$ and the pair $(\epsilon_1,\epsilon_2)\in (\F_q^*)^2$.
According to Lemma 5.1 and Lemma 5.2, the sequence of partial
quotients for $\alpha$
satisfies $(I)$ if and only if there exists a sequence
$(\delta_n)_{n\geq 0}$ in $\F_q^*$ satisfying $(II)_0$, $(D_2)$ and
$(D_3)$, where the sequence $(\lambda_n)_{n\geq 1}$ is based on
$(\delta_n)_{n\geq 0}$ by $(LD)$. 
Given the value for $\delta_0$ in Lemma 5.1, the existence of this sequence requires condition $(II)$
of Theorem B. According to Lemma 5.3, this sequence does exist if and only if there exists a sequence
$(\gamma_n)_{n\geq 1}$ satisfying $(D)$, $(\Gamma_1)$ and $(\Gamma_2)$. Now distinguish two 
cases :either $\epsilon_2^r\delta_l-4k^2\theta_k\epsilon_1^r=0$ or
$\epsilon_2^r\delta_l-4k^2\theta_k\epsilon_1^r\neq 0$. In the first case, which is case $(III)_1$ of Theorem B, 
by $(\Gamma_1)$ and according to the previous value for $\delta_0$, we have $\gamma_1=0$ and 
also, by $(\Gamma_2)$, $\gamma_n=0$ for $n\geq 2$. In the second case, which is case $(III)_2$ of
 Theorem B, again by $(\Gamma_1)$ and $(\Gamma_2)$, the sequence $(\gamma_n)_{n\geq 1}$ 
is in $\F_q^*$ and consequently, using Lemma 5.4, this sequence can be described by the
 formulas $(\Gamma)$ of Theorem B. In both cases, the sequence $(\delta_n)_{n\geq 0}$
is described recursively from $\delta_1,\dots,\delta_l$ and by $(D)$
from the sequence $(\gamma_n)_{n\geq 1}$, identically zero or not. So the proof of the theorem is complete. 
\vskip 0.2 cm
\noindent {\emph{ Proof of Corollary C:}}
 \par First, due to the degrees of the polynomial coefficients of
 equation $(E)$, we observe that if this equation has a root $\alpha$ in
 $\F(27)$ then we must have $\vert \alpha \vert=\vert T \vert$. Now,
 with Proposition A, we consider the continued fraction of type
 $(3,1,1)$, in $\F(27)$
 defined by $\lambda_1=1$ and the pair $(-u^6,u^3)\in
 (\F_{27}^*)^2$. So we have
 $\alpha^3=-u^6(T^2-1)\alpha_2+u^3T$, where $\alpha_2=1/(\alpha-T)$.
 Hence this continued fraction satisfies equation $(E)$. This proves
 that $(E)$ has no other root in $\F(27)$ (and consequently this root
 is algebraic over $\F_{27}(T)$ of degree four). Now we need to prove
 that the expansion for $\alpha$ is perfect. Here we have $k=1$ and
 $l=1$, consequently $\theta_1=1$ and $f(n)=3n-1$ for $n\geq 1$. From Lemma 5.1, we also have
 $\delta_0=u^{-3}$ and $\delta_1=-\lambda_1^3+\epsilon_2=-1+u^3=u^4\in
 \F_{27}^*$. So $(II)$ holds. We compute now $\gamma_1$. We have
 $\gamma_1^3=(4\theta_1\epsilon_1^3\delta_1^{-3}-\epsilon_2^3)\theta_1\delta_1^{-3}=u^3\neq 0$. 
We are in case $(III_2)$ if the sequence $(\gamma_n)_{n\geq 1}$ can be defined. 
First we compute $C_0$. We have $C_0=-\delta_1^{-1}(\gamma_1\epsilon_1^3)(\delta_1\gamma_1)^{-3}=1$.
 Applying the formulas in $(\Gamma)$ with the triplet $(h_0,h_1,h_2)$ in $(\F_3(X))^3$ stated 
in $(H)$, we obtain the recursive definition given in the corollary for $(\gamma_n)_{n\geq 1}$.
 As $\gamma_1=u$ has degree 3 over $\F_3$, by induction all the terms have the same degree 
over $\F_3$ and therefore this sequence $(\gamma_n)_{n\geq 1}$ is well defined. 
Applying the formulas $(D)$ with the triplet $(g_0,g_1,g_2)$ in $(\F_3(X))^3$ stated in $(G)$, 
we obtain the recursive definition for the sequence $(\delta_n)_{n\geq 1}$ from $(\gamma_n)_{n\geq 1}$
 as stated in the corollary. Finally, applying the formulas $(LD)$,
 the sequence $(\lambda_n)_{n\geq 1}$ satisfies the recursive definition 
indicated in the corollary. This completes the proof.

\end{document}